\documentclass[12pt,twoside]{article} 
\usepackage{amsfonts,amssymb,amsmath}
\usepackage{url}

\date{} 
\newtheorem{thrm}{\quad Theorem}[section]
\newtheorem{lem}[thrm]{\quad Lemma}
\newtheorem{prop}[thrm]{\quad Proposition}
\newtheorem{cor}[thrm]{\quad Corollary}

\input{tcilatex}

\begin{document} 

\centerline{\bf  International Mathematical Forum, Vol. x, 200x, no. xx, xxx - xxx} 

\centerline{} 

\centerline{} 

\centerline {\Large{\bf A Systematic Study of Frame Sequence Operators}} 

\centerline{} 

\centerline{\Large{\bf and their Pseudoinverses}} 

\centerline{} 

\centerline{\bf {P. Balazs}} 

\centerline{} 

\centerline{Acoustics Research Institute} 

\centerline{Austrian Academy of Sciences} 

\centerline{Vienna, Austria}

\centerline{} 

\centerline{\bf {M. A. El-Gebeily}} 

\centerline{} 

\centerline{Department of Mathematical Sciences} 

\centerline{King Fahd University of Petroleum and Minerals} 

\centerline{Dhahran, Saudi Arabia} 

\begin{abstract}
In this note we investigate the operators associated with frame sequences in
a Hilbert space $H$, i.e., the synthesis operator $T:\ell ^{2}\left( 
\mathbb{N}
\right) \rightarrow H$ , the analysis operator $T^{\ast }:H\rightarrow $ $%
\ell ^{2}\left( 
\mathbb{N}
\right) $ and the associated frame operator $S=TT^{\ast }$ as operators
 defined on (or to) the whole space rather than on subspaces.
Furthermore, the projection $P$ onto the range of $T$, the projection $Q$
onto the range of $T^{\ast }$ and the Gram matrix $G=T^{\ast }T$ are
investigated. For all these operators, we investigate their pseudoinverses,
how they interact with each other, as well as possible classification of
frame sequences with them. For a tight frame sequence, we show that some of
these operators are connected in a simple way.
\end{abstract}



\section{Introduction}

\bigskip Frame sequences are the natural generalization of frames \cite%
{Casaz1}. In many situations, for example, when constructing a frame
multi-resolution analysis (see e.g., \cite{BenedetoLi}, \cite{Christbk}), we
start with a frame sequence in a Hilbert space $H$ and then define the
initial approximation space $V$ as the closure of the span of the frame
sequence. In the literature, several operators are associated with frame
sequences and these spaces, namely, the projection operator $P:H\rightarrow H
$ onto $V,$ the inclusion operator $\iota _{V}:V\rightarrow H$ , the
analysis operator $\mathcal{U}:V\rightarrow $ $\ell ^{2}\left( 
\mathbb{N}
\right) $, the synthesis operator $\mathcal{T}:\ell ^{2}\left( 
\mathbb{N}
\right) \rightarrow V$ and the frame operator $\mathcal{S}:V\rightarrow V$
(the definition of all these operators will be given in the next section).
In the literature, frame sequences are analyzed mostly using the concrete
representation of these operators. On the other hand, analyzing frame
sequences from a pure operator theoretic point of view can offer a deeper
insight into the structure of such sequences. 
This note will take 
this approach. Almost exclusively, all the
proofs provided in this paper use these operators and operator theoretic
principles.

An alternative way of looking at the above operators is to extend their
definitions with the help of the first two operators to the whole space $H.$
This way we are always working with the base spaces $H$ and $\ell ^{2}\left( 
\mathbb{N}
\right) $ and we do not worry about the subspace $V$ and its image in $\ell
^{2}\left( 
\mathbb{N}
\right) .$ The extended operators become: the synthesis operator $T:\ell
^{2}\left( 
\mathbb{N}
\right) \rightarrow H$, the analysis operator $U:H\rightarrow \ell
^{2}\left( 
\mathbb{N}
\right) $ and the frame operator $S:H\rightarrow H.$ We develop
relationships between these extended versions. In fact, we start by defining
these extended versions rather than the "classical" restricted operators.
Consequently no ambiguity arises regarding notions such as inverses of
operators and pseudo inverses. While the proofs of the relationships are
straightforward in the most part, they are nontrivial in the sense that we
use the current state of knowledge to derive them. 
Also, a form of duality between statements
about the synthesis and the analysis operators emerges throughout the
presentation. For all involved operators, we investigate how they interact
with each other, as well as possible classification of frame sequences with
them. For a tight frame sequence, we will show that some of these operators
are connected in a simple way.

Of course, most of the work in the literature on frame sequences is related
to this work. We mention in particular the references \cite{Chist99}, \cite%
{Chistetal00}, \cite{Teolis93} which have more direct bearing on this note.

As a preliminary lemma we list here the properties of the pseudo-inverse or
the Moore Penrose inverse of a bounded operator with closed range that are
most important to us (see, for example, Appendix A.7 in \cite{Christbk}).

\begin{lem}
\label{lem:pinv}Let $H_{1},H_{2}$ be Hilbert spaces, and suppose that $%
U:H_{1}\rightarrow H_{2}$ is a bounded operator with closed range $R_{U}$.
Then there exists a bounded operator $U^{\dag }:H_{2}\rightarrow H_{1}$ such
that 
\[
UU^{\dag }f=f,\ \forall f\in R_{U}, 
\]

with $\ker _{U^{\dagger }}=\overline{R}_{U}^{\bot }$ and $\overline{R}%
_{U^{\dagger }}=\ker _{U}^{\bot }$. This operator is uniquely determined by
these properties.

Furthermore, $U^{\dag }$ has the following properties.

\begin{enumerate}
\item $UU^{\dag }$ is the orthogonal projection of $H_{2}$ onto $R_{U}$.

\item $U^{\dag }U$ is the orthogonal projection of $H_{1}$ onto $R_{U^{\dag
}}.$

\item $U^{\ast }$ has closed range and $\left( U^{\ast }\right) ^{\dag
}=\left( U^{\dag }\right) ^{\ast }.$

\item $U^{\dagger }UU^{\dagger }=U^{\dagger }$. \label{eq:aux0}

\item On $\overline{R}_U$ we have $U^\dagger = U^* \left( U U^* \right)^{-1}$%
.
\end{enumerate}
\end{lem}

\section{Frame Sequences and Their Pseudoinverses}

In this section $H$ denotes a general Hilbert space and $\ell ^{2}\left( 
\mathbb{N}
\right) $ denotes the space of absolutely square summable sequences of
complex numbers. We will denote elements of $\ell ^{2}\left( 
\mathbb{N}
\right) $ by lower case letters such as $c,d,\cdots $ etc. and, when we want
to explicitly use the terms of the sequences $c,d,\cdots ,$ we will use
Greek letters such as $\left( \zeta _{k}\right) _{k=1}^{\infty },\left( \eta
_{k}\right) _{k=1}^{\infty },\cdots $ etc. When no confusion arises we will
write these sequences as $\left( \zeta _{k}\right) ,\left( \eta _{k}\right)
,\cdots $ etc. We denote by $\left\{ \epsilon _{k}\right\} _{k=1}^{\infty }$
the sequence of standard basis elements in $\ell ^{2}\left( 
\mathbb{N}
\right) .$

Suppose $\left\{ f_{k}\right\} _{k=1}^{\infty }$ is a sequence in $H.$ With $%
\left\{ f_{k}\right\} _{k=1}^{\infty }$ we associate three, possibly
unbounded \cite{antino1}, operators: the synthesis operator $T:\ell
^{2}\left( 
\mathbb{N}
\right) \rightarrow H$ defined by´
$
Tc=\sum_{k=1}^{\infty }\zeta _{k}f_{k}, 
$, %
the analysis operator $U:H\rightarrow \ell ^{2}\left( 
\mathbb{N}
\right) $ defined by 
$
Uf=\left( \left\langle f,f_{k}\right\rangle \right) 
$ 
and the \emph{frame operator} $S$ defined by 
$
Sf=TUf=\sum_{k=1}^{\infty }\left\langle f,f_{k}\right\rangle f_{k} 
$, %
whenever the right hand sides of these definitions exist. Observe that $T$ \
is densely defined as its domain $D_{T}$ contains all finite sequences
(sequences which, eventually, consist of zeros) in $\ell ^{2}\left( 
\mathbb{N}
\right) .$ This implies that $T$ \ has a well defined adjoint $T^{\ast
}:\ell ^{2}\left( 
\mathbb{N}
\right) \rightarrow H$, which is a closed operator (see \cite{Rudin}). We
also have

\begin{description}
\item[(A)] \bigskip\ $\ker _{T^{\ast }}$ is closed.

\item[(B)] $\ker _{T^{\ast }}=\left( R_{T}\right) ^{\bot }=\left( \overline{R%
}_{T}\right) ^{\bot }.$

\item[(C)] $H=\overline{R}_{T}\oplus \left( \overline{R}_{T}\right) ^{\bot }=%
\overline{R}_{T}\oplus $ $\ker _{T^{\ast }}.$
\end{description}

It follows that the orthogonal projection $P$ of $H$ onto $\overline{R}_{T}=$
$\left( \ker _{T^{\ast }}\right) ^{\bot }$ is always well defined. The
following lemma and its corollary are straightforward.

\begin{lem}
$\limfunc{span}\left\{ f_{k}\right\} _{k=1}^{\infty }\subseteq
R_{T}\subseteq \overline{\limfunc{span}}\left\{ f_{k}\right\} _{k=1}^{\infty
}.$
\end{lem}

\begin{cor}
\label{cor:closed range} If $T$ has closed range, then $R_{T}=\overline{%
\limfunc{span}}\left\{ f_{k}\right\} _{k=1}^{\infty }.$
\end{cor}

Recall that $\left\{ f_{k}\right\} _{k=1}^{\infty }$ is a frame sequence in $%
H$ if there are positive constants $A,B$ such that%
\begin{equation}
A\left\Vert f\right\Vert ^{2}\leq \sum_{k=1}^{\infty }\left\vert
\left\langle f,f_{k}\right\rangle \right\vert ^{2}\leq B\left\Vert
f\right\Vert ^{2}\ \forall f\in \limfunc{span}\left\{ f_{k}\right\}
_{k=1}^{\infty }.  \label{eq:framineq}
\end{equation}%
The frame sequence is tight if $A=B.$ By the definition above, it becomes
clear that it is a frame for $V=\overline{\limfunc{span}}\{f_{k}\}$. The
restricted versions of the above operators are defined in the same way with
the only difference that they work from or to $V.$ We let $\iota
_{V}:V\rightarrow H$ be the inclusion operator $\iota _{V}(f)=f$ , $\mathcal{%
U}:V\rightarrow $ $\ell ^{2}\left( 
\mathbb{N}
\right) $ the analysis operator, $\mathcal{T}:\ell ^{2}\left( 
\mathbb{N}
\right) \rightarrow V$ the synthesis and $\mathcal{S}:V\rightarrow V$ the
frame operator. We have following basic relationships between these
operators, which are straightforward to show

\begin{prop}
If $\left\{ f_{k}\right\} _{k=1}^{\infty }$ is a frame sequence in $H$, then
the following properties hold.

\begin{enumerate}
\item $T=\iota _{V}\mathcal{T}$.

\item $\overline{R}_{T}=\overline{R}_{\mathcal{T}}=V$ and $P$ is the
projection on $V$.

\item $U=\mathcal{U}P$.

\item $R_{U}=R_{\mathcal{U}}.$

\item $S=TU=\iota _{V}\mathcal{T}\mathcal{U}P=\iota _{V}\mathcal{S}P$.
\end{enumerate}
\end{prop}

\begin{proof}
(iii) \& (Iv) : For $h_{2}\in V^{\bot },$ we clearly have $U(h_{2})=\left(
\left\langle h_{2},f_{k}\right\rangle \right) =0$. Every $h\in H$ can be
uniquely be described as $h=h_{1}+h_{2}$ with $h_{1}\in V$ and $h_{2}\in
V^{\bot }$, therefore, $U(h)=U(h_{1})+U(h_{2})=U(h_{1})=\mathcal{U}Pf.$ 

All
the others proofs are straightforward.
\end{proof}

As a frame sequence is a frame for its closed span, $\mathcal{U}$ and $%
\mathcal{T}$ are bounded, $\mathcal{T=U}^{\ast }$ and $\mathcal{U=T}^{\ast }$%
. Also

\begin{cor}
If $\left\{ f_{k}\right\} _{k=1}^{\infty }$ is a frame sequence, then

\begin{enumerate}
\item the analysis operator $U$ is bounded,

\item the synthesis operator $T$ is bounded,

\item $T=U^{\ast }$ and $U=T^{\ast }.$
\end{enumerate}
\end{cor}

\begin{proof}
(i) and (ii): Since $\mathcal{T}$, $\mathcal{U}$, $P$ and $\iota _{V}$ are
bounded, the boundedness of $U$ and $T$ follows directly from the above
relations. 

(iii): $U^{\ast }=\left( \mathcal{U}P\right) ^{\ast }=P^{\ast }%
\mathcal{U}^{\ast }=\iota _{V}\mathcal{T}=T$. Repeat the argument for $%
T^{\ast }$.
\end{proof}

As $\{f_{k}\}$ is a frame for $V$, there is a sequence $\{\tilde{f}%
_{k}\}\subseteq V$, which is the canonical dual frame for $V$, $\tilde{f}%
_{k}=\mathcal{S}^{-1}f_{k}$. $\{\tilde{f}_{k}\}$ is again a frame sequence
in $H$ with closed span $V$ and the bounds $\tilde{A}=\frac{1}{B}$ and $%
\tilde{B}=\frac{1}{A}$. Let $\tilde{T},\tilde{U},\tilde{S},\tilde{\mathcal{T}%
},\tilde{\mathcal{U}}$ and $\tilde{\mathcal{S}}$ be the corresponding
operators associated with this frame sequence. Therefore we have $\tilde{%
\mathcal{T}}=\mathcal{S}^{-1}\mathcal{T}$ , $\ \tilde{\mathcal{U}}=\mathcal{%
US}^{-1}$, $\tilde{\mathcal{S}}=\mathcal{S}^{-1}$ and $\mathcal{T}\tilde{%
\mathcal{U}}=\tilde{\mathcal{T}}\mathcal{U}=id_{V}$. The following corollary
can be easily shown.

\begin{cor}
If $\left\{ f_{k}\right\} _{k=1}^{\infty }$ is a frame sequence in $H$ and $%
\{\tilde{f}_{k}\}$ is its dual sequence, then the following properties hold.

\begin{enumerate}
\item $\tilde{T}=i_{V}\tilde{\mathcal{T}}$ and $\tilde{U}=\tilde{\mathcal{U}}%
P.$

\item $\tilde{S}=\iota _{V}\mathcal{S}^{-1}P.$

\item $T\tilde{U}=\iota _{V}P = \tilde{T} U $.
\end{enumerate}
\end{cor}

It follows from Property (3) above that the projection $P$ on $V$ as a
function from $H$ into $H$ is $T\tilde{U}.$ This is well known \cite%
{Christbk}. Also, denoting by $Q$ the orthogonal projection of $\ell
^{2}\left( 
\mathbb{N}
\right) $ onto $\left( \ker _{T}\right) ^{\bot }=\overline{R}_{T^{\ast },}$
it is straightforward to show that $Q$ is the Gram matrix $G=U\tilde{T}=%
\mathcal{U}\tilde{\mathcal{T}}$ .

We have the following characterizations of frame sequences.

\begin{thrm}
\label{thm:mainsec2}The following are equivalent:

\begin{enumerate}
\item $\left\{ f_{k}\right\} _{k=1}^{\infty }$ is a frame sequence in $H$
with bounds $A,B.$

\item There exist positive constants $A,B$ such that for every $f\in H,$ 
$$
A\left\Vert Pf\right\Vert ^{2}\leq \left\Vert T^{\ast }f\right\Vert ^{2}\leq
B\left\Vert Pf\right\Vert ^{2}.  \label{eq:ineq1}
$$

\item There exist positive constants $A,B$ such that for every $c\in \ell
^{2}\left( 
\mathbb{N}
\right) ,$%
$$
A\left\Vert Qc\right\Vert ^{2}\leq \left\Vert Tc\right\Vert ^{2}\leq
B\left\Vert Qc\right\Vert ^{2}.  \label{eq:ineq2}
$$
\end{enumerate}
\end{thrm}

\begin{proof}
Assume (i) holds. As $T^{\ast }=U=\mathcal{U}P$, (ii) is equivalent to the
definition of frame sequences and is therefore true.

Assume (ii) holds. Then a similar inequality holds for the dual frame in $V$,
i.e., 
$$\tilde{A}\left\Vert Pf\right\Vert ^{2}\leq \left\Vert \tilde{U}%
f\right\Vert ^{2}\leq \tilde{B}\left\Vert Pf\right\Vert ^{2}, \ \mathrm{or} \ \frac{1}{B}%
\left\Vert Pf\right\Vert ^{2}\leq \left\Vert \tilde{U}f\right\Vert ^{2}\leq 
\frac{1}{A}\left\Vert Pf\right\Vert ^{2}.$$
Now choose $c\in \ell ^{2}\left( 
\mathbb{N}
\right) $ and set $f=Tc$. Then $\frac{1}{B}\left\Vert PTc\right\Vert
^{2}\leq \left\Vert \tilde{U}Tc\right\Vert ^{2}\leq \frac{1}{A}\left\Vert
PTc\right\Vert ^{2},$ which implies that 
$$\frac{1}{B}\left\Vert Tc\right\Vert
^{2}\leq \left\Vert Qc\right\Vert ^{2}\leq \frac{1}{A}\left\Vert
Tc\right\Vert ^{2}, \ \mathrm{ or } \ A\left\Vert Qc\right\Vert ^{2}\leq \left\Vert
Tc\right\Vert ^{2}\leq B\left\Vert Qc\right\Vert ^{2}.$$

Assume (iii) holds. Since $\left\Vert Qc\right\Vert \leq $ $\left\Vert
c\right\Vert ,$ it follows that $T$ is continuous. Let 
$c\in \ker _{T}^{\bot }$. Then 
\[
A\left\Vert c\right\Vert ^{2}\leq \left\Vert Tc\right\Vert ^{2}\leq
B\left\Vert c\right\Vert ^{2}. 
\]%
Therefore $T|_{\ker _{T}^{\bot }}$ is bounded, injective and has closed
range \cite{conw1,xxlphd1}. As $R_{T}=R_{T|_{\ker _{T}^{\bot }}}$, $T$ is
bounded and has closed range. Using \cite{Christppr} this is equivalent to $%
\{f_{k}\}$ forming a frame sequence.
\end{proof}

It follows easily from Theorem \ref{thm:mainsec2} that $\left\{
f_{k}\right\} _{k=1}^{\infty }$ is a frame for $H$ if and only if $T$ is
surjective (i.e., $P=I$) and that $\left\{ f_{k}\right\} _{k=1}^{\infty }$
is a Riesz Basis if and only if $T^{\ast }$ is onto (i.e., $Q=I$). Theorem %
\ref{thm:mainsec2} can also be restated in a number of other ways; the
following one uses the optimal bounds for the inequalities. 

\begin{cor}
\label{cor:5}The following are equivalent:

\begin{enumerate}
\item $\left\{ f_{k}\right\} _{k=1}^{\infty }$ is a frame sequence in $H.$

\item $T^{\ast }$ is continuous, has a closed range, and%
\[
\left\Vert T^{\dag }\right\Vert ^{-1}\left\Vert Pf\right\Vert \leq
\left\Vert T^{\ast }f\right\Vert \leq \left\Vert T\right\Vert \left\Vert
Pf\right\Vert \ \forall f\in H. 
\]

\item $T$ is continuous, has a closed range, and%
\[
\left\Vert T^{\dag }\right\Vert ^{-1}\left\Vert Qc\right\Vert \leq
\left\Vert Tc\right\Vert \leq \left\Vert T\right\Vert \left\Vert
Qc\right\Vert \ \forall c\in \ell ^{2}\left( 
\mathbb{N}
\right) . 
\]
\end{enumerate}
\end{cor}

For a frame sequence the frame operator $S = T T^*$ is bounded and self
adjoint. Furthermore, since $R_{T^{\ast }}$ is closed, 
\begin{eqnarray*}
R_{S} &=&TT^{\ast }H=T\left( T^{\ast }H+\ker _{T}\right) \\
&=&T\left( T^{\ast }H+R_{T^{\ast }}^{\bot }\right) =T\left( R_{T^{\ast
}}+R_{T^{\ast }}^{\bot }\right) \\
&=&T\ell ^{2}\left( 
\mathbb{N}
\right) =R_{T}.
\end{eqnarray*}%
It follows that $S$ has closed range and, hence, a continuous Moore-Penrose
pseudo-inverse $S^{\dagger }.$ We list some properties of $S^{\dagger }$
next.

\begin{lem}
\label{lem:S properties}The operator $S^{\dagger }:H\rightarrow H$ is the
same as the operator $\tilde{S}=i_{V}\mathcal{S}^{-1}P$ and therefore is
self-adjoint and has the following properties:

\begin{enumerate}
\item $SS^{\dagger }= S^{\dagger }S = P$ $.$

\item $S^{\dagger }\left( I-P\right) =0.$

\item $S^{\dagger }P=PS^{\dagger }=S^{\dagger }.$
\end{enumerate}
\end{lem}

\begin{proof}
$\tilde{T}$ has the same range as $T$ as mentioned above, $\overline{R}_{%
\tilde{T}}=V$. Therefore, $\overline{R}_{\tilde{S}}=\overline{R}_{S}$. Also, 
$S\tilde{S}=\iota _{V}\mathcal{S}P\iota _{V}\tilde{\mathcal{S}}P=\iota _{V}%
\mathcal{S}\tilde{\mathcal{S}}P=\iota _{V}\mathcal{SS}^{-1}P=\iota _{V}P.$ 
Furthermore, $\ker {\tilde{S}}=\left\{ f:\tilde{S}f=0\right\} =\left\{
f:\iota _{V}\mathcal{S}Pf=0\right\} =\left\{ f:Pf=0\right\} =V^{\bot },$ as $%
\iota _{V}$ and $\mathcal{S}$ are injective. Repeat the same argument for
the frame sequence $\left\{ \widetilde{f}_{k}\right\} $ with the roles of $S$
and $\tilde{S}$ switched and use Lemma \ref{lem:pinv} to arrive at $%
S^{\dagger }=\tilde{S}$.

(i) By Property (i) of Lemma \ref{lem:pinv}, $SS^{\dagger }$ is the
orthogonal projection onto $R_{S}=R_{T}.$ Therefore, $SS^{\dagger }=P$ $.$
Switch the roles of $S$ and $\tilde{S}$ to show the second part.

(ii) $S^{\dagger }\left( I-P\right) =S^{\dagger }-S^{\dagger
}P=S^{\dagger }-S^{\dagger }SS^{\dagger }=S^{\dagger }-S^{\dagger }=0.$

(iii) The equality 
$ S^{\dagger }P=S^{\dagger }  \label{eq:aux1}
$ follows form (ii). To show that $PS^{\dagger }=S^{\dagger }$ observe first
that $S^{\dagger }P:H\rightarrow R_{S}.$ Hence, by \ref{eq:aux1}, Lemma \ref%
{lem:pinv} and (1.), $S^{\dagger }=S^{\dagger }P=PS^{\dagger }P=PS^{\dagger
}SS^{\dagger }=PS^{\dagger }.$ 
\end{proof}

Because of the frame property on $V,$ among all sequences $c\in \ell
^{2}\left( 
\mathbb{N}
\right) $ which synthesize an $f\in H,$ the sequence $c_{0}=\left(
\left\langle f,S^{\dag }f_{k}\right\rangle \right) $ is the one with the
minimum norm.

Similarly, among all elements $f\in H$ which analyze to a $c\in \ell
^{2}\left( 
\mathbb{N}
\right) ,$ the element $f_{0}=S^{\dagger }Tc=\sum_{k=1}^{\infty }\zeta
_{k}S^{\dag }f_{k}$ is the one with the minimum norm. We have $\left\Vert
f\right\Vert ^{2}=\left\Vert f_{0}\right\Vert ^{2}+\left\Vert
f-f_{0}\right\Vert ^{2}.$

Proposition 5.3.5 in \cite{Christbk} can now be restated in terms of $S^{\dag
}$ which is defined on all of $H$ instead of $\mathcal{S}^{-1}$ which is
defined only on $V$.

\begin{cor}
Let $\left\{ f_{k}\right\} _{k=1}^{\infty }$ be a frame sequence in $H$ .
For any $f\in H,$%
\[
Pf=\sum_{k=1}^{\infty }\left\langle f,S^{\dag }f_{k}\right\rangle f_{k}. 
\]
\end{cor}

\begin{prop}
The pseudo-inverse of $T$ is $\tilde{U}$, $T^{\dagger }=\tilde{U}$. The
pseudo-inverse of $U$ is $\tilde{T}$, $U^{\dagger }=\tilde{T}$.
Consequently, we have the following properties

\begin{enumerate}
\item $T^{\dag }=T^{\ast }S^{\dag }$ and $T^{\dag }= S^{\dag } T^{\ast }$

\item $\left( T^{\dagger }\right) ^{\ast }T^{\dagger }=S^{\dagger }.$

\item $\left( T^{\dagger }\right) ^{\ast }=S^{\dagger }T.$

\item $\left\Vert T\right\Vert ^{2}=\left\Vert S\right\Vert .$

\item $\left\Vert T^{\dagger }\right\Vert ^{2}=$ $\left\Vert S^{\dagger
}\right\Vert .$
\end{enumerate}
\end{prop}

\begin{proof}
Clearly $T \tilde{U} = \iota_V P$ and $ker_{\tilde{U}} = V^\bot$. Again by
switching the roles of $T$ and $\tilde{T}$ we arrive with Lemma \ref%
{lem:pinv} at $T^\dagger = \tilde{U}$. Use an analog argument for $U$ and $%
\tilde{U}$.

(i) $\tilde{U}=\iota _{V}\widetilde{\mathcal{U}}P=\iota _{V}\mathcal{\ 
\widetilde{\mathcal{U}}\left( \tilde{T}U\right) }P=\iota _{V}\mathcal{\
\left( \widetilde{\mathcal{U}}\tilde{T}\right) U}P=\iota _{V}\mathcal{\ 
\tilde{S}U}P=\iota _{V}\mathcal{\ \tilde{S}}P\iota _{V}\mathcal{U}P=\tilde{S}%
U$, so $T^{\dagger }=S^{\dagger }T^{\ast }$. Since $P=SS^{\dag }=S\tilde{S}%
=\left( TT^{\ast }\right) \left( \tilde{T}\tilde{T}^{\ast }\right) =T\cdot %
\left[ T^{\ast }\left( TT^{\ast }\right) ^{\dag }\right] ,$ it follows from
the uniqueness of the Moor-Penrose pseudo-inverse that $T^{\dag }=T^{\ast
}S^{\dag }.$

(ii) This follows immediately from $T^{\dagger }=\tilde{U}$

(iii) follows from (i) by taking conjugates.

(iv): We have for every $f\in H,\qquad \left\Vert T^{\ast }f\right\Vert
^{2}=\left\langle Sf,f\right\rangle \leq \left\Vert Sf\right\Vert \left\Vert
f\right\Vert \leq \left\Vert S\right\Vert \left\Vert f\right\Vert ^{2}.$On
the other hand, $\left\Vert S\right\Vert =\left\Vert TT^{\ast }\right\Vert
\leq \left\Vert T\right\Vert ^{2}.$

(v): Observe first that, for every $f\in H,$ 
\begin{equation}
\left\Vert T^{\dag }f\right\Vert ^{2}=\left\langle f,S^{\dagger
}f\right\rangle .  \label{eq:aux3}
\end{equation}%
This can be seen as follows: 
\begin{eqnarray*}
\left\Vert T^{\dag }f\right\Vert ^{2} &=&\left\langle T^{\ast }S^{\dagger
}f,T^{\ast }S^{\dagger }f\right\rangle =\left\langle TT^{\ast }S^{\dagger
}f,S^{\dagger }f\right\rangle  \\
&=&\left\langle Pf,S^{\dagger }f\right\rangle =\left\langle f,PS^{\dagger
}f\right\rangle =\left\langle f,S^{\dagger }f\right\rangle .
\end{eqnarray*}%
Hence, for every $f\in H,$ $\left\Vert T^{\dag }f\right\Vert ^{2}\leq
\left\Vert S^{\dagger }\right\Vert \left\Vert f\right\Vert ^{2}.$It follows
that $\left\Vert T^{\dag }\right\Vert ^{2}\leq \left\Vert S^{\dagger
}\right\Vert .$On the other hand, since $S^{\dagger }$ is self adjoint, 
\[
\left\Vert S^{\dagger }\right\Vert =\sup_{\left\Vert f\right\Vert
=1}\left\langle S^{\dagger }f,f\right\rangle =\sup_{\left\Vert f\right\Vert
=1}\left\Vert T^{\dagger }f\right\Vert ^{2}\leq \left\Vert T^{\dagger
}\right\Vert ^{2}.
\]%
Therefore, $\left\Vert S^{\dagger }\right\Vert =\left\Vert T^{\dagger
}\right\Vert ^{2}.$
\end{proof}

Let $\left\{ f_{k}\right\} _{k=1}^{\infty }$ be a frame sequence in $H.$
Define the Gram matrix $G:\ell ^{2}\left( 
\mathbb{N}
\right) \rightarrow $ $\ell ^{2}\left( 
\mathbb{N}
\right) $ by $G=UT=T^{\ast }T.$Alternatively, $G=\mathcal{U}Pi_{V}\mathcal{T}%
=\mathcal{UT}.$ More explicitly, $Gc=\sum_{l=1}^{\infty } \left( \sum_{k=1}^{\infty }c_{k}\left\langle
f_{k},f_{j}\right\rangle \right) \epsilon _{j}$ and $\left\langle Gc,c\right\rangle
=\sum_{j,k=1}^{\infty }\left\langle c_{k}f_{k},c_{j}f_{j}\right\rangle .$It
immediately follows from the definitions that $T^{\ast }S=GT^{\ast }$ and $%
ST=TG.$ Clearly $G$ is self adjoint. It is well known \cite{xxlphd1,Christbk}
that $G$ is a bijective bounded operator from $R_{T^{\ast }}$ onto $R_{T}$
with bounded inverse if and only if $\left\{ f_{k}\right\} $ is a frame
sequence. In particular $R_{G}=R_{T}$ and $\ker _{G}=\ker _{T}=R_{T^{\ast
}}^{\bot }$. It follows that $G$ has a closed range and, hence, a continuous
Moor-Penrose pseudo-inverse $G^{\dagger }.$ The Gram matrix is the
projection onto $\left( \ker _{T}\right) ^{\bot }=\overline{R}_{T^{\ast }}.$

We list some properties of $G^{\dagger }$ next. For that let us denote by $%
\tilde{G}$ the Gram matrix corresponding to the dual frame $\left\{ 
\widetilde{f}_{k}\right\} .$The proof is analogous to the one of Lemma \ref%
{lem:S properties} with appropriate adjustments.

\begin{lem}
\label{lem:St properties} The operator $G^{\dagger }:\ell ^{2}\left( 
\mathbb{N}
\right) \rightarrow \ell ^{2}\left( 
\mathbb{N}
\right) $ is the same as $\tilde{G}$. It is therefore self-adjoint and has
the following properties:

\begin{enumerate}
\item $G G^{\dagger }=Q = G^{\dagger } G$ $.$

\item $G^{\dagger }\left( I-Q\right) =0.$

\item $G^{\dagger }Q=QG^{\dagger }=G^{\dagger }.$
\end{enumerate}
\end{lem}

The following corollary is the same as \cite{Christbk}, Proposition 5.3.6.

\begin{cor}
Let $\left\{ f_{k}\right\} _{k=1}^{\infty }$ be a frame sequence in $H$ .
For any $c\in \ell ^{2}\left( 
\mathbb{N}
\right) ,Qc=\sum_{k=1}^{\infty }\left\langle c,G^{\dag }T^{\ast
}f_{k}\right\rangle \epsilon _{k}.$
\end{cor}

We may also write $Qc=\sum_{k,j=1}^{\infty }\left\langle c,G^{\dag }\epsilon
_{j}\right\rangle \left\langle f_{j},f_{k}\right\rangle \epsilon _{k}.$

\begin{lem}
We have the following properties

\begin{enumerate}
\item $\left( T^{\ast }\right) ^{\dag }=T G^{\dag }.$

\item $T^{\dagger }\left( T^{\dagger }\right) ^{\ast }=G^{\dagger }.$

\item $T^{\dagger }=G^{\dagger }T^{\ast }.$

\item $\left\Vert T\right\Vert ^{2}=\left\Vert G \right\Vert .$

\item $\left\Vert T^{\dagger }\right\Vert ^{2}=$ $\left\Vert G^{\dagger
}\right\Vert .$
\end{enumerate}
\end{lem}

\begin{proof}
(i): Since $Q=T^{\ast }T G^{\dag }=$ $T^{\ast }\left( T^{\ast }\right)
^{\dag },$ it follows from the uniqueness of the Moor-Penrose pseudo-inverse
that $\left( T^{\ast }\right) ^{\dag }=T G^{\dag }. $

(ii) follows immediately from (1.) by multiplying on the left by $T^{\dagger
}$ and using Property 3 of Lemma \ref{lem:pinv}.

(iii) follows from (1.) by taking adjoints.

(iv): We have for every $c\in \ell ^{2}\left( 
\mathbb{N}
\right) ,$\ $\left\Vert Tc\right\Vert ^{2}=\left\langle Gc,c\right\rangle
\leq \left\Vert G\right\Vert \left\Vert c\right\Vert ^{2}.$ Thus $\left\Vert
T\right\Vert ^{2}\leq \left\Vert G\right\Vert .$On the other hand, $%
\left\Vert G\right\Vert =\left\Vert T^{\ast }T\right\Vert \leq \left\Vert
T\right\Vert ^{2},$ which yeilds $\left\Vert G\right\Vert =\left\Vert
T\right\Vert ^{2}.$

(v): As $T^{\dagger }=\tilde{T}$ and $G^{\dagger }=\tilde{G}$, (5) is the
same as (4) for the dual frame.
\end{proof}

\begin{cor}
$\left\Vert G \right\Vert =\left\Vert S\right\Vert $and $\left\Vert
G^{\dagger }\right\Vert =\left\Vert S^{\dagger }\right\Vert .$
\end{cor}

Theorem \ref{thm:mainsec2} (or rather, Corollary \ref{cor:5}) can also be
reformulated as

\begin{thrm}
The following are equivalent:

\begin{enumerate}
\item $\left\{ f_{k}\right\} _{k=1}^{\infty }$ is a frame sequence in $H.$

\item $S$ is continuous, has closed range and%
\[
\left\Vert S^{\dag }\right\Vert ^{-1}\left\Vert Pf\right\Vert ^{2}\leq
\left\langle Sf,f\right\rangle \leq \left\Vert S\right\Vert \left\Vert
Pf\right\Vert ^{2}\ \forall f\in H. 
\]

\item $G$ is continuous, has closed range and%
\[
\left\Vert S^{\dagger }\right\Vert ^{-1}\left\Vert Qc\right\Vert ^{2}\leq
\left\langle G c,c\right\rangle \leq \left\Vert S\right\Vert \left\Vert
Qc\right\Vert ^{2}\ \forall c\in \ell ^{2}\left( 
\mathbb{N}
\right) . 
\]
\end{enumerate}
\end{thrm}

\begin{lem}
Let $\left\{ f_{k}\right\} _{k=1}^{\infty }$ be a tight frame sequence in $%
H. $ Then

\begin{enumerate}
\item $S=AP,\ G=AQ,$

\item $\ S^{\dag }=\frac{1}{A}P,\ G^{\dagger }=\frac{1}{A}Q.$
\end{enumerate}
\end{lem}

\begin{proof}
We prove the statements for $G$ only. Since the frame is tight, $%
\left\langle Gc,c\right\rangle =\left\Vert Tc\right\Vert ^{2}=A\left\Vert
Qc\right\Vert ^{2}=A\left\langle Qc,Qc\right\rangle .$
\end{proof}

An analog of the polarization identity can be easily proved:%
\[
\left\langle Gc,d\right\rangle =%
\begin{array}{c}
\frac{1}{4}\left( \left\langle \ G\left( c+d\right) ,c+d\right\rangle
-\left\langle \ G\left( c-d\right) ,c-d\right\rangle \right. + \\ 
\left. +i\left\langle \ G\left( c+id\right) ,c+id\right\rangle
-i\left\langle \ G\left( c-id\right) ,c-id\right\rangle \right) ,%
\end{array}%
\]%
which yields%
\begin{eqnarray*}
\left\langle Gc,d\right\rangle  &=&%
\begin{array}{c}
\frac{1}{4}A\left( \ \left\Vert Q\left( c+d\right) \right\Vert
^{2}-\left\Vert Q\left( c-d\right) \right\Vert ^{2}\right.+   \\ 
\left. +i\left\Vert Q\left( c+id\right) \right\Vert ^{2}-i\left\Vert Q\left(
c-id\right) \right\Vert ^{2}\right) 
\end{array}
\\
&=&A\left\langle Qc,Qd\right\rangle =A\left\langle Qc,d\right\rangle .
\end{eqnarray*}%
Since $d$ is arbitrary, $Gc=AQc.$ Furthermore, $Qc=\ G^{\dagger
}Gc=G^{\dagger }AQc=AG^{\dagger }c.$ This gives $G^{\dagger }c=\frac{1}{A}Qc.
$

\vspace{1cm}

{\bf ACKNOWLEDGEMENTS.} \\
The authors would like to thank Pete Casazza for creating the contact
between them, as well as giving comments on a first version of
this paper.

The work of the first author was partly supported by the European Union's Human Potential Program, under contract HPRN-CT-2002-00285 (HASSIP). He would like to thank the hospitality of the LATP, CMI, Marseille, France and FYMA, UCL, Louvain-la-Neuve, Belgium.

{\bf Received: August 28, 2007}

\end{document}